\documentclass{article}
\usepackage{graphicx}
\usepackage{nips07submit_e,times}
\usepackage{subfigure}
\usepackage{amsmath}
\usepackage{amssymb}
\usepackage{float}

\title{Upper-Bounding Proof Length with the Busy Beaver}

\author{Gustavo Lacerda\\
\texttt{gusl@optimizelife.com}}

\usepackage{amsmath, amssymb, graphics, latexsym, makeidx, graphicx}
\begin{document}

\maketitle

%% \mainmatter

%% \paper[Short Theorems must have ``Short'' Proofs]{An Information-Theoretic Upper Bound on the Length of the Shortest Proof}{Gustavo Lacerda}{Machine Learning Department\\
%% Carnegie Mellon University\\
%% Address\\
%% Pittsburgh, PA, 15213 USA\\
%% %%\texttt{gusl@cs.cmu.edu}}{}

\begin{abstract}
Consider a short theorem, i.e. one that can be written down using just a few symbols. Can its shortest proof be arbitrarily long? We answer this question in the negative. Inspired by arguments by Calude et al (1999) and Chaitin (1984) that construct an upper bound on the first counterexample of a $\Pi_1$ sentence as a function of the sentence's length, we present a similar argument about proof length for arbitrary statements. As with the above, our bound is uncomputable, since it uses a Busy Beaver oracle. Unlike the above, our result is not restricted to any complexity class.
Finally, we combine the above search procedures into an automatic (albeit uncomputable) procedure for discovering G\"{o}del sentences.
\end{abstract}

\section{Introduction}

Suppose you have a hypothesis of the form $\forall x . \phi(x)$, where $x \in \mathbb{N^+}$ and $\phi$ is a decidable predicate. Statements of this form are known as a $\Pi_1$ statements, and this class includes famous problems such as the Goldbach conjecture (``all even numbers are sums of two primes'').

%Collatz problem? ``does x's sequence end in 1?'' is perhaps only semi-decidable
% but interesting to bound the number of iterations with BB

It is commonly taught that checking examples never suffices to establish the truth of a hypothesis, {\it no matter how many examples one has checked}, and that this is why a {\it proof} is needed. However, Chaitin (1984) has shown an upper bound on the smallest counterexample (should it exist), based on the length used to encode the statement. This bound uses the Busy Beaver function, as shown in the next section, in which we reproduce Chaitin's argument.

This paper presents an analogous result about proofs. It is likewise commonly held that no matter how hard you've tried and failed to prove a result, it is always possible that the proof is out there, and you just haven't found it yet. This note shows an upper bound on proof length, which essentially says that if no proof can be found within a certain finite set of proof attempts, then no proof exists.

We conclude by combining the above two results into a brief discussion of G\"{o}del statements and when to introduce new axioms.

It is important now to remark that a Busy Beaver oracle is equivalent to a Halting oracle. As a result, the bounds discussed here are uncomputable, and are likely to be unknown for any given hypothesis that one might encode (and possibly unknowable by all axiom systems in current use). Therefore, these results might only begin to be useful if we ever have estimates of BB for large enough integers.

\section{Bounding the first counterexample}

The following argument is due to Chaitin (1984) and Calude, J\"{u}rgensen \& Legg (1999).

Let the program $Ch$ be a checker: given as inputs a predicate $\phi$ over the positive natural numbers, and a positive natural number $x$, $Ch(\phi, x)$ always halts with output telling us whether or not $\phi(x)$ is true.

\begin{figure}
\begin{verbatim}
function Ch(phi,x)
  return phi(x)
end function
\end{verbatim}
\end{figure}
Let $P(\phi)$ be the program that returns the smallest $x$ on which $Ch(\phi,x)$ returns false. It works by counting up until it reaches a number for which $Ch$ returns false.

\begin{verbatim}
function phi(x)
  ...
end

function P(phi)
  x <- 0
  while (phi(x)==TRUE)
    x <- x+1
  end while
  return x
end function
\end{verbatim}

Thus, for any hypothesis $\forall x \phi(x)$,  $P(\phi)$ either gives us the smallest counterexample (in case the hypothesis is false) or it never halts (in case the hypothesis is true).

%We have exhibited a program for computing the smallest counterexample for any $\Pi_1$ statement $\forall x \phi(x)$ over the naturals, namely $P(\phi)$. 

%%Note that we haven't written the entire program here, since it requires implementing the function \verb+phi+.

The size of this program is $length(\ulcorner \phi \urcorner) + length(\ulcorner P \urcorner)$ \footnote{We are omitting the log term for brevity. See Li \& Vitanyi for an explanation of why this is needed.}, where $\ulcorner f \urcorner$ refers to the implementation of function $f$. Here, $length(\ulcorner P \urcorner)$ isn't very large, since it fits in 5 lines of code.

 Therefore, should a counterexample to $\forall x \phi(x)$ exist, the Kolmogorov complexity of the smallest counterexample is $\le length(\ulcorner \phi \urcorner) + length(\ulcorner P \urcorner)$.

Therefore, in order to decide the truth of the universal, we only need to check the numbers whose Kolmogorov Complexity $\le length(\ulcorner \phi \urcorner) + length(\ulcorner P \urcorner)$. The Busy Beaver function gives us an upper bound \cite{Cha84}: we only need to check numbers up to $BB(length(\ulcorner \phi \urcorner) + length(\ulcorner P \urcorner))$: if no counterexample is found, none exist.

That is, the smallest counterexample is $\le BB(length(\ulcorner \phi \urcorner) + length(\ulcorner P \urcorner))$.

\section{Bounding the size of the smallest proof}

Now we turn to our original result:

Let the program $TP$ be a theorem pump (a breadth-first search over all proofs): given a sentence $T$, it will search until it finds a proof of $T$. $TP(T)$ will either halt with the shortest proof (if there is a proof), or never halt (if there isn't).

The size of $TP(T)$ is $length(\ulcorner T \urcorner)+length(\ulcorner TP \urcorner)$. Thus the KC of the shortest proof is $\le length(\ulcorner T \urcorner)+length(\ulcorner TP \urcorner)$. Given this KC, the Busy Beaver function gives us an upper bound on the length of the shortest proof: $BB(length(\ulcorner T \urcorner) + length(\ulcorner TP \urcorner))$. So we only need to check proofs up to that size: if none is found, none exist.

\section{Conclusion}

In systems that are sound and complete, one could use either approach to decide the truth for any given $\Pi_1$ sentence. However, in interesting theories (namely, those that express PA), it follows from G\"{o}del's 2nd incompleteness theorem that there will be sentences that are true in PA but not provable in the theory.

Together, the above two procedures (counterexample search and proof search) can be combined into a procedure for discovering $\Pi_1$ G\"{o}del sentences for any given axiom set $T$, i.e. sentences that are unprovable in $T$ (since no proof was found within our bound) but true in PA (since no counterexample was found within the Chaitin bound).

\subsection*{Acknowledgements}
The author thanks Shane Legg for comments.

\section*{References}

G. J. Chaitin (1984) - Computing the Busy Beaver function. In {\it Information, Randomness \& Incompleteness, pages 74-76}

C. Calude, H. J\"{u}rgensen and S. Legg (1999) - Solving Problems with Finite Test Sets. In {\it Finite versus Infinite: Contributions to an Eternal Dilemma, pages 39-52, Springer-Verlag, London}

M. Li, P. Vitanyi (1997) - An Introduction to Kolmogorov Complexity and Its Applications. {\it New York: Springer-Verlag}.

Solomon Feferman - Transfinite Recursive Progressions of Axiomatic Theories {\it The Journal of Symbolic Logic, Vol. 27, No. 3 (Sep., 1962), pp. 259-316}

%Does incompleteness theorem suggest?

%Feferman (1960) suggests extending theories by adding  G\"{o}del sentences as new axioms, and points out that G\"{o}del's 2nd incompleteness theorem implies that this is an infinite hierarchy.

%%%%%%% This is where your paper is included. %%%%%%
%%\include{template_bib}

%%\backmatter

\end{document}